\documentclass[a4paper,11pt]{amsart}
\usepackage{amsmath,amssymb,amsfonts,latexsym}

\newtheorem{theorem}{Theorem}[section]

\newtheorem{corollary}[theorem]{Corollary}

\input{tcilatex}

\begin{document}
\title[\textbf{New} \textbf{Generalization of Eulerian polynomials}]{\textbf{%
New generalization of Eulerian polynomials and their applications }}
\author[\textbf{S. Araci}]{\textbf{Serkan Araci}}
\address{\textbf{University of Gaziantep, Faculty of Science and Arts,
Department of Mathematics, 27310 Gaziantep, TURKEY}}
\email{\textbf{mtsrkn@hotmail.com}}
\author[\textbf{M. Acikgoz}]{\textbf{Mehmet Acikgoz}}
\address{\textbf{University of Gaziantep, Faculty of Science and Arts,
Department of Mathematics, 27310 Gaziantep, TURKEY}}
\email{\textbf{acikgoz@gantep.edu.tr}}
\author[\textbf{E. \c{S}en}]{\textbf{Erdo\u{g}an \c{S}en}}
\address{\textbf{Department of Mathematics, Faculty of Science and Letters,
Nam\i k Kemal University, 59030 Tekirda\u{g}, TURKEY}}
\email{\textbf{erdogan.math@gmail.com}}

\begin{abstract}
In the present paper, we introduce Eulerian polynomials with parameters $a$
and $b$ and give the definition of them. By using the definition of
generating function for our polynomials, we derive some new identities in
Analytic Numbers Theory. Also, we give relations between Eulerian
polynomials with parameters $a$ and $b$, Bernstein polynomials,
Poly-logarithm functions, Bernoulli and Euler numbers. Moreover, we see that
our polynomials at $a=-1$ are related to Euler-Zeta function at negative
inetegers. Finally, we get Witt's formula for new generalization of Eulerian
polynomials which we express in this paper.

\vspace{2mm}\noindent \textsc{2010 Mathematics Subject Classification.}
Primary 05A10, 11B65; Secondary 11B68, 11B73.

\vspace{2mm}

\end{abstract}

\maketitle




\section{\textbf{Introduction}}


The Bernoulli numbers and polynomials, Euler numbers and polynomials,
Genocchi numbers and polynomials, Stirling numbers of the second kind,
Bernstein polynomials and Eulerian polynomials possess many interesting
properties not only in Complex Analysis, and Analytic Numbers Theory but
also in Mathematical physics, and $p$-adic analysis. These polynomials have
been studied by many mathematicians for a long time (for details about this
subject, see [1-28]).

Eulerian polynomial sequence $\left\{ \mathcal{A}_{n}\left( x\right)
\right\} _{n\geq 0}$ is given by the following summation:%
\begin{equation}
\sum_{l=0}^{\infty }l^{n}x^{l}=\frac{\mathcal{A}_{n}\left( x\right) }{\left(
1-x\right) ^{n+1}},\text{ \ }\left\vert x\right\vert <1.  \label{Equation 1}
\end{equation}

It is well-known that the Eulerian polynomial, $\mathcal{A}_{n}\left(
x\right) $, of degree $n$ can be introduced as%
\begin{equation}
\mathcal{A}_{n}\left( x\right) =\sum_{k=1}^{n}\mathcal{A}\left( n,k\right)
x^{k},\text{ \ }\mathcal{A}_{0}\left( z\right) =1  \label{Equation 2}
\end{equation}

where $\mathcal{A}\left( n,k\right) $ are called the Eulerian numbers that
can be computed by using%
\begin{equation}
\mathcal{A}\left( n,k\right) =\sum_{j=0}^{k}\binom{n+1}{j}\left( -1\right)
^{j}\left( k-j\right) ^{n},\text{ \ }1\leq k\leq n,  \label{Equation 3}
\end{equation}

where $\mathcal{A}\left( n,0\right) =1$. Eulerian polynomials, $\mathcal{A}%
_{n}\left( x\right) $, are also given by means of the following exponential
generating function:%
\begin{equation}
e^{\mathcal{A}\left( x\right) t}=\sum_{n=0}^{\infty }\mathcal{A}_{n}\left(
x\right) \frac{t^{n}}{n!}=\frac{1-x}{e^{t\left( 1-x\right) }-x}
\label{Equation 4}
\end{equation}

where $\mathcal{A}^{n}\left( x\right) :=\mathcal{A}_{n}\left( x\right) $,
symbolically. Eulerian polynomials can be found via the following recurrence
relation:%
\begin{equation}
\left( \mathcal{A}\left( t\right) +\left( t-1\right) \right) ^{n}-t\mathcal{A%
}_{n}\left( t\right) =\left\{ 
\begin{array}{cc}
1-t, & \text{if }n=0 \\ 
0, & \text{if }n\neq 0\text{,}%
\end{array}%
\right.  \label{Equation 5}
\end{equation}

(for details, see \cite{Kim 3}, \cite{Kim 4}, \cite{Foata}, \cite{Araci 2}
and \cite{Araci 4}).

Now also, we give the definition of Eulerian fraction, $\alpha _{n}\left(
x\right) $, can be expressed as%
\begin{equation}
\alpha _{n}\left( x\right) :=\frac{\mathcal{A}_{n}\left( x\right) }{\left(
1-x\right) ^{n+1}}\text{.}  \label{Equation 6}
\end{equation}

We want to note that Eulerian fraction is very useful in the study of the
Eulerian numbers, Eulerian polynomials, Euler function and its
generalization, Jordan function in Number Theory (for details, see \cite%
{Birkhoff}).

In 2010, Acikgoz and Araci firstly introduced the generating function of
Bernstein polynomials, as follows:%
\begin{equation}
\sum_{n=k}^{\infty }B_{k,n}\left( x\right) \frac{t^{n}}{n!}=\frac{\left(
tx\right) ^{k}}{k!}e^{t\left( 1-x\right) },\text{ }t\in 
\mathbb{C}
\text{,}  \label{Equation 7}
\end{equation}%
where $B_{k,n}\left( x\right) $ are called Bernstein polynomials, which are
defined by%
\begin{equation}
B_{k,n}\left( x\right) =\binom{n}{k}x^{k}\left( 1-x\right) ^{n-k},\text{ }%
0\leq x\leq 1,  \label{Equation 8}
\end{equation}

(for details on this subject, see \cite{acikgoz 1}).

The Poly-logarithms can be defined by the series:%
\begin{equation}
Li_{n}\left( z\right) =\sum_{k=1}^{\infty }\frac{z^{k}}{k^{n}}
\label{Equation 9}
\end{equation}

for $n\geq 0$, and $\left\vert z\right\vert <1$. We easily see, for $n=0$%
\begin{equation*}
Li_{0}\left( z\right) =\frac{z}{1-z}.
\end{equation*}

Also, Poly-logarithms can be given by the integral representation, as
follows:%
\begin{equation*}
Li_{n}\left( z\right) =\int_{0}^{z}\frac{Li_{n}\left( z\right) }{z}dz
\end{equation*}%
in $%
\mathbb{C}
\backslash \left[ 1,\infty \right) .$ We note that $Li_{1}\left( z\right)
=-\log \left( 1-z\right) $ is the usual logarithm (see \cite{Lewin}).

In \cite{Luo}, \cite{Luo 1} and \cite{Luo 2}, Luo $et$ $al.$ defined the
generalization of the Bernoulli and Euler polynomials with parameters $a,b,c$
as follows:%
\begin{eqnarray}
\frac{tc^{xt}}{b^{t}-a^{t}} &=&\sum_{n=0}^{\infty }\frac{B_{n}\left(
x;a,b,c\right) }{n!}t^{n},\text{ }\left( \left\vert t\log \frac{b}{a}%
\right\vert <2\pi \right) \text{,}  \label{Equation 10} \\
\frac{2c^{xt}}{b^{t}+a^{t}} &=&\sum_{n=0}^{\infty }\frac{E_{n}\left(
x;a,b,c\right) }{n!}t^{n},\text{ }\left( \left\vert t\log \frac{b}{a}%
\right\vert <\pi \right) \text{.}  \label{Equation 11}
\end{eqnarray}

So that, obviously,%
\begin{equation}
B_{n}\left( x;1,e,e\right) :=B_{n}\left( x\right) \text{ and }E_{n}\left(
x;1,e,e\right) :=E_{n}\left( x\right) .  \label{Equation 12}
\end{equation}

Here $B_{n}\left( x\right) $ and $E_{n}\left( x\right) $ are the classical
Bernoulli polynomials and the classical Euler polynomials, respectively.

Next, for the classical Bernoulli numbers, $B_{n}$ and the classical Euler
numbers, $E_{n}$ we have%
\begin{equation}
B_{n}\left( 0\right) :=B_{n}\text{ and }E_{n}\left( 0\right) :=E_{n}.
\label{Equation 13}
\end{equation}

By the same motivation of all the above generalizations, we consider, in
this paper, the generalization of Eulerian polynomials and derive some new
theoretical properties for them. Also, we show that our polynomials are
related to poly-logarithm function, the Bernstein polynomials, Bernoulli
numbers, Euler numbers, Genocchi numbers, Euler-Zeta function and Stirling
numbers of the second kind. Finally, we get Witt's formula for new
generalization of Eulerian polynomials which seems to be interesting for
further work in $p$-adic analysis.

\section{\textbf{On the new generalization of Eulerian polynomials}}

In this section, we start by giving the following definition of new
generalization of Eulerian polynomials.

\begin{definition}
Let $b\in 
\mathbb{R}
^{+}$(positive real numbers) and $a\in 
\mathbb{C}
$(field of complex numbers), then we define the following:%
\begin{equation}
e^{t\mathcal{A}\left( a,b\right) }=\sum_{n=0}^{\infty }\mathcal{A}_{n}\left(
a,b\right) \frac{t^{n}}{n!}=\frac{1-a}{b^{t\left( 1-a\right) }-a}
\label{Equation 14}
\end{equation}%
where $\mathcal{A}_{n}\left( a,b\right) $ are called the generalization of
Eulerian polynomials (or Eulerian polynomials with parameters $a$ and $b$).
Also, $\mathcal{A}^{n}\left( a,b\right) :=\mathcal{A}_{n}\left( a,b\right) $%
, symbolically. 
\end{definition}

So that, obviously,%
\begin{equation*}
\mathcal{A}_{n}\left( x,e\right) :=\mathcal{A}_{n}\left( x\right) \text{.}
\end{equation*}

By (\ref{Equation 14}), we have the following recurrence relation for the
Eulerian polynomials with parameters $a$ and $b$:%
\begin{equation*}
e^{t\mathcal{A}\left( a,b\right) }=\frac{1-a}{e^{t\left( 1-a\right) \ln b}-a}%
.
\end{equation*}

By applying combinatorial techniques to the above equality, then we easily
derive the following theorem:

\begin{theorem}
The following recurrence relation holds: 
\begin{equation}
\left[ \mathcal{A}_{n}\left( a,b\right) +\left( 1-a\right) \ln b\right]
^{n}-a\mathcal{A}_{n}\left( a,b\right) =\left( 1-a\right) \delta _{n,0}
\label{Equation 15}
\end{equation}%
where $\delta _{n,0}$ is the Kronecker's symbol.
\end{theorem}

We now consider for $n>0$ in (\ref{Equation 15}), becomes%
\begin{equation}
\mathcal{A}_{n}\left( a,b\right) =\frac{1}{a-1}\sum_{k=0}^{n-1}\binom{n}{k}%
\mathcal{A}_{k}\left( a,b\right) \left( 1-a\right) ^{n-k}\left( \ln b\right)
^{n-k}.  \label{Equation 16}
\end{equation}

We want to note that taking $a=x$ and $b=e$ in (\ref{Equation 16}) reduces
to 
\begin{equation}
\mathcal{A}_{n}\left( x\right) =\frac{1}{x-1}\sum_{k=0}^{n-1}\binom{n}{k}%
\mathcal{A}_{k}\left( x\right) \left( 1-x\right) ^{n-k}  \label{Equation 17}
\end{equation}%
(see \cite{Kim 3} and \cite{Foata}). We see that (\ref{Equation 17}) is
proportional with Bernstein polynomials which we state in the following
theorem:

\begin{theorem}
The following identity%
\begin{equation*}
\mathcal{A}_{n}\left( x\right) =\sum_{k=0}^{n-1}\frac{\mathcal{A}_{k}\left(
x\right) B_{k,n}\left( x\right) }{x^{k+1}-x^{k}}
\end{equation*}%
is true.
\end{theorem}

Let us now consider $\lim_{t\rightarrow 0}\frac{d^{k}}{dt^{k}}$ in (\ref%
{Equation 14}), then we readily arrive at the following theorem.

\begin{theorem}
Let $b\in 
\mathbb{R}
^{+}$ and $a\in 
\mathbb{C}
$, then we have%
\begin{equation}
\mathcal{A}_{k}\left( a,b\right) =\lim_{t\rightarrow 0}\left[ \frac{d^{k}}{%
dt^{k}}\left( \frac{1-a}{b^{t\left( 1-a\right) }-a}\right) \right] \text{.}
\label{Equation 18}
\end{equation}
\end{theorem}

By (\ref{Equation 18}), we easily conclude the following corollary.

\begin{corollary}
The following Cauchy-type integral holds true:%
\begin{equation*}
\frac{1}{1-a}\mathcal{A}_{k}\left( a,b\right) =\frac{k!}{2\pi i}\int_{C}%
\frac{t^{-k-1}}{b^{t\left( 1-a\right) }-a}dt
\end{equation*}%
where $C$ is a loop which starts at $-\infty ,$ encircles the origin once in
the positive direction, and the returns $-\infty $.
\end{corollary}

By (\ref{Equation 14}), we discover the following:%
\begin{eqnarray*}
\sum_{n=0}^{\infty }\mathcal{A}_{n}\left( a^{2},b^{2}\right) \frac{t^{n}}{n!}
&=&\left[ \frac{1-a}{b^{t\left( 1+a\right) \left( 1-a\right) }-a}\right] %
\left[ \frac{1+a}{b^{t\left( 1-a\right) \left( 1+a\right) }-a}\right] \\
&=&\left[ \sum_{n=0}^{\infty }\left( 1+a\right) ^{n}\mathcal{A}_{n}\left(
a,b\right) \frac{t^{n}}{n!}\right] \left[ \sum_{n=0}^{\infty }\left(
1-a\right) ^{n}\mathcal{A}_{n}\left( -a,b\right) \frac{t^{n}}{n!}\right] .
\end{eqnarray*}

By using Cauchy product on the above equality, then we get the following
theorem.

\begin{theorem}
The following equality%
\begin{equation}
\mathcal{A}_{n}\left( a^{2},b^{2}\right) =\sum_{k=0}^{n}\binom{n}{k}\left(
1+a\right) ^{k}\mathcal{A}_{k}\left( a,b\right) \mathcal{A}_{n-k}\left(
-a,b\right) \left( 1-a\right) ^{n-k}  \label{Equation 19}
\end{equation}%
is true.
\end{theorem}

After the basic operations in (\ref{Equation 19}), we discover the following
corollary.

\begin{corollary}
The following property holds:%
\begin{equation*}
\mathcal{A}_{n}\left( a^{2},b^{2}\right) =\sum_{k=0}^{n}\left( 1+\frac{1}{a}%
\right) ^{k}B_{k,n}\left( a\right) \mathcal{A}_{k}\left( a,b\right) \mathcal{%
A}_{n-k}\left( -a,b\right) \text{.}
\end{equation*}
\end{corollary}

Now also, we consider geometric series in (\ref{Equation 14}), then we
compute as follows:%
\begin{eqnarray*}
\sum_{n=0}^{\infty }\mathcal{A}_{n}\left( a,b\right) \frac{t^{n}}{n!} &=&%
\frac{1-a}{e^{t\left( 1-a\right) \ln b}-a}=\frac{1-a^{-1}}{%
1-a^{-1}e^{t\left( 1-a\right) \ln b}} \\
&=&\left( 1-\frac{1}{a}\right) \sum_{j=0}^{\infty }a^{-j}e^{jt\left(
1-a\right) \ln b} \\
&=&\left( 1-\frac{1}{a}\right) \sum_{j=0}^{\infty }a^{-j}\sum_{n=0}^{\infty
}j^{n}\left( 1-a\right) ^{n}\left( \ln b\right) ^{n}\frac{t^{n}}{n!} \\
&=&\sum_{n=0}^{\infty }\left[ \left( 1-\frac{1}{a}\right) \sum_{j=0}^{\infty
}a^{-j}j^{n}\left( 1-a\right) ^{n}\left( \ln b\right) ^{n}\right] \frac{t^{n}%
}{n!}.
\end{eqnarray*}

By comparing the coefficients of $\frac{t^{n}}{n!}$ on the above equation,
then we readily derive the following theorem.

\begin{theorem}
The following%
\begin{equation*}
\left( \frac{1}{a-1}\right) ^{n}\mathcal{A}_{n}\left( a,b\right) =\left[ 
\frac{\left( \ln b\right) ^{n}}{a}-\left( \ln b\right) ^{n}\right]
\sum_{j=1}^{\infty }\frac{a^{-j}}{j^{-n}}
\end{equation*}%
is true.
\end{theorem}

The above theorem is related to Poly-logarithm function, as follows:%
\begin{equation}
\left( \frac{1}{a-1}\right) ^{n}\mathcal{A}_{n}\left( a,b\right) =\left[ 
\frac{\left( \ln b\right) ^{n}}{a}-\left( \ln b\right) ^{n}\right]
Li_{-n}\left( a^{-1}\right) \text{.}  \label{Equation a}
\end{equation}

In \cite{Lewin}, it is well-known that%
\begin{equation}
Li_{-n}\left( x\right) =\left( x\frac{d}{dx}\right) ^{n}\frac{x}{1-x}%
=\sum_{k=0}^{n}k!S\left( n+1,k+1\right) \left( \frac{x}{1-x}\right) ^{k+1}
\label{Equation aa}
\end{equation}

where $S\left( n,k\right) $ are the Stirling numbers of the second kind. By (%
\ref{Equation a}) and (\ref{Equation aa}), we have the following interesting
theorem.

\begin{theorem}
The following holds true:%
\begin{equation*}
a\mathcal{A}_{n}\left( a,b\right) =-\left( \ln b\right)
^{n}\sum_{k=0}^{n}k!S\left( n+1,k+1\right) \left( \frac{1}{a-1}\right) ^{k-n}%
\text{.}
\end{equation*}
\end{theorem}

\section{\textbf{Further Remarks}}

Now, we consider (\ref{Equation 14}) for evaluating at $a=-1$, as follows:%
\begin{equation}
\sum_{n=0}^{\infty }\mathcal{A}_{n}\left( -1,b\right) \frac{t^{n}}{n!}=\frac{%
2}{b^{2t}+1}  \label{Equation 20}
\end{equation}%
where $\mathcal{A}_{n}\left( -1,b\right) $ are called Eulerian polynomials
with parameter $b$.

By (\ref{Equation 20}), we derive the following equality in complex plane:%
\begin{equation*}
\sum_{n=0}^{\infty }i^{n}\mathcal{A}_{n}\left( -1,b\right) \frac{t^{n}}{n!}=%
\frac{2}{b^{2it}+1}=\frac{2}{e^{2it\ln b}+1}.
\end{equation*}

From this, we discover the following:%
\begin{equation}
\sum_{n=0}^{\infty }i^{n}\mathcal{A}_{n}\left( -1,b\right) \frac{t^{n}}{n!}%
=\sum_{n=0}^{\infty }E_{n}2^{n}i^{n}\left( \ln b\right) \frac{t^{n}}{n!}
\label{Equation 21}
\end{equation}%
where $E_{n}$ are $n$-th Euler numbers which are defined by the following
exponential generating function:%
\begin{equation}
\sum_{n=0}^{\infty }E_{n}\frac{t^{n}}{n!}=\frac{2}{e^{t}+1}\text{, }%
\left\vert t\right\vert <\pi \text{.}  \label{Equation 22}
\end{equation}

By (\ref{Equation 21}) and (\ref{Equation 22}), we have the following
theorem.

\begin{theorem}
Let $n\in 
\mathbb{N}
$ (field of natural numbers) and $b\in 
\mathbb{C}
$, then we get%
\begin{equation*}
\mathcal{A}_{n}\left( -1,b\right) =2^{n}E_{n}\left( \ln b\right) ^{n}\text{.}
\end{equation*}
\end{theorem}

We now give the definition of Bernoulli numbers for sequel of this paper via
the following exponential generating function:%
\begin{equation}
\sum_{n=0}^{\infty }B_{n}\frac{t^{n}}{n!}=\frac{t}{e^{t}-1}\text{, }%
\left\vert t\right\vert <2\pi \text{.}  \label{Equation 23}
\end{equation}

By using (\ref{Equation 20}) and (\ref{Equation 23}), we see that%
\begin{equation*}
\sum_{n=0}^{\infty }\mathcal{A}_{n}\left( -1,b\right) \frac{t^{n}}{n!}=\frac{%
2}{e^{2t\ln b}+1}=\frac{1}{t}\left[ \frac{2t}{e^{2t\ln b}-1}-\frac{4t}{%
e^{4t\ln b}-1}\right] \text{.}
\end{equation*}

So from above%
\begin{equation*}
\sum_{n=0}^{\infty }\mathcal{A}_{n}\left( -1,b\right) \frac{t^{n}}{n!}%
=\sum_{n=0}^{\infty }\left[ 2^{n}\left( \ln b\right) ^{n}B_{n}-4^{n}\left(
\ln b\right) ^{n}B_{n}\right] \frac{t^{n-1}}{n!}\text{.}
\end{equation*}

By comparing the coefficients of $t^{n}$ on the above equation, then we can
state the following theorem.

\begin{theorem}
The following identity%
\begin{equation*}
\mathcal{A}_{n}\left( -1,b\right) =\frac{2^{n+1}\left( \ln b\right)
^{n+1}\left( 1-2^{n+1}\right) B_{n+1}}{n+1}
\end{equation*}%
holds true.
\end{theorem}

By (\ref{Equation 20}), we obtain the following:%
\begin{equation*}
\sum_{n=0}^{\infty }\mathcal{A}_{n}\left( -1,b\right) \frac{t^{n}}{n!}=\frac{%
1}{t}\sum_{n=0}^{\infty }2^{n}\left( \ln b\right) ^{n}G_{n}\frac{t^{n}}{n!}%
=\sum_{n=0}^{\infty }2^{n}\left( \ln b\right) ^{n}G_{n}\frac{t^{n-1}}{n!}%
\text{.}
\end{equation*}

That is, we reach the following theorem.

\begin{theorem}
The following holds true:%
\begin{equation*}
\mathcal{A}_{n}\left( -1,b\right) =\frac{2^{n+1}\left( \ln b\right)
^{n+1}G_{n+1}}{n+1}
\end{equation*}%
where $G_{n}$ are the familiar Genocchi numbers which is defined by%
\begin{equation*}
\sum_{n=0}^{\infty }G_{n}\frac{t^{n}}{n!}=\frac{2t}{e^{t}+1}\text{.}
\end{equation*}
\end{theorem}

We reconsider (\ref{Equation 20}) and using definition of geometric series,
then we compute as follows:%
\begin{eqnarray*}
\sum_{n=0}^{\infty }\mathcal{A}_{n}\left( -1,b\right) \frac{\left( \frac{t}{2%
}\right) ^{n}}{n!} &=&2\sum_{j=0}^{\infty }\left( -1\right) ^{j}e^{jt\ln b}
\\
&=&\sum_{n=0}^{\infty }\left( 2\left( \ln b\right) ^{n}\sum_{j=0}^{\infty
}\left( -1\right) ^{j}j^{n}\right) \frac{t^{n}}{n!}\text{.}
\end{eqnarray*}

Therefore, we obtain the following theorem

\begin{theorem}
For $n>0$, then we have 
\begin{equation}
\mathcal{A}_{n}\left( -1,b\right) =2^{n+1}\left( \ln b\right)
^{n}\sum_{j=1}^{\infty }\left( -1\right) ^{j}j^{n}\text{.}
\label{Equation 24}
\end{equation}
\end{theorem}

As is well known, Euler-zeta function is defined by 
\begin{equation}
\zeta _{E}\left( s\right) =2\sum_{j=1}^{\infty }\frac{\left( -1\right) ^{j}}{%
j^{s}},\text{ }s\in 
\mathbb{C}
\text{.}  \label{Equation 25}
\end{equation}

From (\ref{Equation 24}) and (\ref{Equation 25}), we obtain the
interpolation function of new generalization of Eulerian polynomials at $%
a=-1 $, as follow:%
\begin{equation}
\mathcal{A}_{n}\left( -1,b\right) =2^{n}\left( \ln b\right) ^{n}\zeta
_{E}\left( -n\right) \text{.}  \label{Equation 26}
\end{equation}

Equation (\ref{Equation 26}) seems to be interpolation function at negative
integers for Eulerian polynomials with parameter $b$.

Let us now consider Witt's formula for our polynomials at $a=-1,$ so we need
the following notations:

Imagine that $p$ be a fixed odd prime number. Throughout this paper, we use
the following notations. By $%
\mathbb{Z}
_{p}$, we denote the ring of $p$-adic rational integers, $%
\mathbb{Q}
$ denotes the field of rational numbers, $%
\mathbb{Q}
_{p}$ denotes the field of $p$-adic rational numbers, and $%
\mathbb{C}
_{p}$ denotes the completion of algebraic closure of $%
\mathbb{Q}
_{p}$. Let $%
\mathbb{N}
$ be the set of natural numbers and $%
\mathbb{N}
^{\ast }=%
\mathbb{N}
\cup \left\{ 0\right\} $.

The $p$-adic absolute value is defined by 
\begin{equation*}
\left\vert p\right\vert _{p}=\frac{1}{p}\text{.}
\end{equation*}

Let $q$ be an indeterminate with $\left\vert q-1\right\vert _{p}<1$.

Let $UD\left( 
\mathbb{Z}
_{p}\right) $ be the space of uniformly differentiable functions on $%
\mathbb{Z}
_{p}$. For a positive integer $d$ with $\left( d,p\right) =1$, let 
\begin{equation*}
X=X_{d}=\lim_{\overleftarrow{n}}%
\mathbb{Z}
/dp^{n}%
\mathbb{Z}
=\underset{a=0}{\overset{dp-1}{\cup }}\left( a+dp%
\mathbb{Z}
_{p}\right) 
\end{equation*}%
with%
\begin{equation*}
a+dp^{n}%
\mathbb{Z}
_{p}=\left\{ x\in X\mid x\equiv a\left( \func{mod}dp^{n}\right) \right\} 
\end{equation*}%
where $a\in 
\mathbb{Z}
$ satisfies the condition $0\leq a<dp^{n}$ and let $\mathcal{\sigma }%
:X\rightarrow 
\mathbb{Z}
_{p}$ be the transformation introduced by the inverse limit of the natural
transformation%
\begin{equation*}
\mathbb{Z}
/dp^{n}%
\mathbb{Z}
\mapsto 
\mathbb{Z}
/p^{n}%
\mathbb{Z}
\text{.}
\end{equation*}

If $f$ is a function on $%
\mathbb{Z}
_{p}$, then we will utilize the same notation to indicate the function $%
f\circ \mathcal{\sigma }$.

For a continuous function $f:X\rightarrow 
\mathbb{C}
_{p},$ the $p$-adic fermionic integral on $%
\mathbb{Z}
_{p}$ is defined by T. Kim in \cite{Kim 1} and \cite{Kim 2}, as follows:

\begin{equation}
I_{-1}\left( f\right) =\int_{X}f\left( \upsilon \right) d\mu _{-1}\left(
\upsilon \right) =\int_{%
\mathbb{Z}
_{p}}f\left( \upsilon \right) d\mu _{-1}\left( \upsilon \right)
=\lim_{n\rightarrow \infty }\sum_{\upsilon =0}^{p^{n}-1}\left( -1\right)
^{\upsilon }f\left( \upsilon \right) \text{.}  \label{Equation 27}
\end{equation}

By (\ref{Equation 27}), it is well-known that%
\begin{equation}
I_{-1}\left( f_{1}\right) +I_{-1}\left( f\right) =2f\left( 0\right) 
\label{Equation 28}
\end{equation}%
where $f_{1}\left( \upsilon \right) :=f\left( \upsilon +1\right) .$
Substituting $f\left( \upsilon \right) =b^{2\upsilon t}$ into (\ref{Equation
28}), we get the following:%
\begin{equation}
\int_{X}e^{2t\upsilon \ln b}d\mu _{-1}\left( \upsilon \right) =\frac{2}{%
b^{2t}+1}=\sum_{n=0}^{\infty }\mathcal{A}_{n}\left( -1,b\right) \frac{t^{n}}{%
n!}.  \label{Equation 29}
\end{equation}

By (\ref{Equation 29}) and using Taylor expansion of $e^{2t\upsilon \ln b}$,
we obtain Witt's formula for our polynomials at $a=-1$, as follows:

\begin{theorem}
The following holds true:%
\begin{equation}
\mathcal{A}_{n}\left( -1,b\right) =\left( \ln b\right)
^{n}2^{n}\int_{X}\upsilon ^{n}d\mu _{-1}\left( \upsilon \right) \text{.}
\label{Equation 30}
\end{equation}
\end{theorem}

Equation (\ref{Equation 30}) seems to be interesting for our further works
in $p$-adic analysis and Analytic numbers theory.

%


\end{document}